\documentclass[a4paper,amsart]{article}
\usepackage{amsmath, amsthm, amssymb}
\usepackage[latin1]{inputenc}
\usepackage[T1]{fontenc}
\usepackage[dvips]{epsfig}
\usepackage{color}

\newtheorem{Proposition}{Proposition}[section]
\newtheorem{Corollary}{Corollary}[section]
\newtheorem{Theorem}{Theorem}[section]
\newtheorem{Remark}{Remark}[section]

\numberwithin{equation}{section}
\def\dvg{{\rm div}}

\def\Rd{{\Bbb R}^d}
\def\Md{{\Bbb M}^{\,d\times d}}
\def\proofbegin{\emph{Proof. }}
\def\proofend{\qedsymbol}
\newcommand\be{\begin{eqnarray*}}
\newcommand\ee{\end{eqnarray*}}
\newcommand\ben{\begin{eqnarray}}
\newcommand\een{\end{eqnarray}}
\newcommand{\comment}[1]{}
\def\NT{\mid\!\mid\!\mid}

\def\dvg{{\rm div}}

\def\IntO{\int\limits_\Omega}

\newcommand{\IntG}[1]{\int\limits_{\Gamma_#1}}
\newcommand{\bup}[1]{\overline{c}_{#1}}
\newcommand{\blow}[1]{\underline{c}_{#1}}
\def\indset{\mathcal{S}}

\begin{document}

\begin{center}
\Large{Two-sided estimates of the solution set for the
reaction-diffusion problem with uncertain data} \\
\vspace{0.7cm}
\large{O. Mali\; and\; S. Repin} \\
\vspace{0.7cm}
\end{center}


\begin{center}
{\sf Dedicated to Prof. R. Glowinski in occasion of his 70th jubilee.}
\end{center}

\begin{abstract}
\textbf{Abstract.} We consider linear reaction--diffusion problems
with mixed Diricl\'et-Neumann-Robin conditions. The diffusion matrix,
reaction coefficient, and the coefficient in the Robin boundary
condition are defined with an uncertainty which allow bounded variations
around some given mean values. A solution to such a problem
cannot be exactly determined (it is a function in the set of
``possible solutions'' formed by generalized solutions related
to possible data). The problem is to find parameters of this
set. In this paper, we show that computable lower and upper
bounds of the diameter (or radius) of the set can be expressed
throughout problem data and parameters that regulate the indeterminacy
range. Our method is based on using a posteriori error majorants
and minorants of the functional type (see \cite{RE,RE1}), which explicitly depend on the
coefficients and allow to obtain the corresponding lower and
upper bounds by solving the respective extremal problems
generated by indeterminacy of coefficients.
\end{abstract}

\section{Introduction} \label{intro}

This paper is concerned with boundary-value problems
for partial differential equations of elliptic type coefficients of
which contain an indeterminacy. Such a situation is quite
typical in real-life problems where parameters of mathematical
model cannot be determined exactly and instead one knows only that the
coefficients belong to a certain set of ``admissible'' data $\Lambda$.
  In view of this fact,
instead of a single exact solution ``$u$'', we have to consider a
``set of solutions'' (we denote it by $\indset( \Lambda )$). As a
result, the  error control problem comes in a more complicated
form in which approximation errors must be evaluated together with
errors arose due to indeterminant data (various approaches that
can be used for such an analysis are exposed in, e.g.,
\cite{HlChBa,MaRe1,MaRe2,RE3}).

In this paper, we
establish explicit relations between the sets $\Lambda$
and $\indset(\Lambda)$ for the reaction-diffusion problem
with mixed Dirichl\'et--Robin boundary conditions
(we call it $\mathcal{P}$) defined by the system
\ben
\label{basic1}
-\dvg( A \nabla u ) + \rho u  & = & f   \quad \textrm{ in } \Omega   \\
\label{basic2}
                           u  & = & 0\quad \textrm{ on } \Gamma_1 \\
\label{basic3}
          n \cdot A \nabla u  & = & F   \quad \textrm{ on } \Gamma_2 \\
\label{basic4}
\alpha u + n \cdot A \nabla u & = & G   \quad \textrm{ on } \Gamma_3 .
\een
Here, $\Omega\in\Rd$ is a bounded and connected domain with Lipschitz
continuous boundary
$\Gamma_1 \cup \Gamma_2 \cup \Gamma_3$ and $f \neq 0$.
We assume that exact  $A$, $\rho$, and $\alpha$ are unknown. Instead,
we know that $A\in \Lambda_A$, $\rho\in\Lambda_\rho$, and $\alpha\in\Lambda_\alpha$,
where
\be
\Lambda_A & := &  \{ A \in L_2(\Omega,\Md) \; | \; A = A_0 + \delta_1 \Psi, \; ||\Psi||_{L_\infty(\Omega,\Md)} \leq 1 \} \\
\Lambda_\rho & := &  \{ \rho \in L_2(\Omega) \; | \; \rho = \rho_0 + \delta_2 \psi_\rho , \; ||\psi_\rho||_{L_\infty(\Omega)} \leq 1 \} \\
\Lambda_\alpha & := &  \{ \alpha \in L_2(\Gamma_3) \; | \; \alpha = \alpha_0 + \delta_3 \psi_\alpha , \; ||\psi_\alpha||_{L_\infty(\Gamma_3)} \leq 1 \} .
\ee
In other words, we assume that the sets of admissible data
are formed  by some (limited) variations of some
known ``mean'' data (which are denoted by subindex 0). The parameters
 $\delta_i$, $i=1,2,3$
represent the magnitude of these variations.
Thus, in the case considered,
\be
\Lambda:=\Lambda_A\times \Lambda_\rho\times \Lambda_\alpha.
\ee
We note that such a presentation of the data arises in
many engineering problems
where data are given in a form \verb+mean+$\pm$\verb+error+. The
solution associated to non-perturbed data $A_0$, $\rho_0$, and $\alpha_0$
 is denoted by $u_0$.

Our goal is to give computable estimates of the radius of $\indset(\Lambda)$,
(we denote this quantity by $r_{\indset}$). The value of $r_{\indset}$ has
a large significance for practical applications because
it shows an accuracy limit defined by the problem statement.
Attempts to find approximate solutions having approximation
errors lesser then $r_{\indset}$ have no practical sense.

The generalized statement of Problem ($\mathcal P$)
 is as follows:
 Find $u \in V_0$ such that
\begin{equation} \label{varform}
  a(u,w) = l(w) \quad \forall w \in V_0 ,
\end{equation}
where space $V_0$ and functionals $a: V_0 \times V_0 \rightarrow R$
and  $l: V_0 \rightarrow R$ are defined by the relations
\be
V_0    & := & \{ w \in H^1(\Omega) \; | \; w_{|_{\Gamma_1}} =  0 \} ,  \\
a(u,w) & := & \IntO  A \nabla u \cdot \nabla w \; dx + \IntO \rho uw \; dx + \IntG{3} \alpha uw \; ds, \\
l(w) & := & \IntO fw \; dx + \IntG{2} Fw \; ds + \IntG{3} Gw \; ds .
\ee
We assume that
\be
\blow{1} |\xi|^2 \leq & A_0 \xi \cdot \xi & \leq \bup{1} |\xi|^2
\quad \forall \xi \in \Rd \quad \textrm{ on } \Omega ,\\
\blow{2} \leq & \rho_0 & \leq \bup{2} \quad \textrm{ on } \Omega, \\
\blow{3} \leq & \alpha_0 & \leq \bup{3} \quad \textrm{ on } \Gamma_3 ,
\ee
where $\blow{i}>0$.
In view of the above-stated conditions, the ``mean'' problem is evidently
elliptic and has a unique solution $u_0$. The condition
\be
0 \leq \delta_i < \blow{i},\qquad i=1,2,3
\ee
guarantees that  the perturbed problem remains elliptic and possesses
a unique solution $u$.

For each $A,\rho,\alpha\in \Lambda$, the corresponding
 problem of $\mathcal{P}(\Lambda)$ is natural to analyze using the (energy) norm
\begin{equation}
\label{1.2}
\NT v \NT_{A,\rho,\alpha}^2 :=
a(v,v)= \IntO A \nabla v \cdot \nabla v \; dx + \IntO \rho v^2 \; dx + \IntG{3} \alpha v^2 \; ds .
\end{equation}
For the sake of simplicity we will also use an abridged notation
$\NT v \NT$ for the norm $\NT v \NT_{A,\rho,\alpha}$.
For the ``mean'' problem,
we use the norm \mbox{$\NT v \NT_{A_0,\rho_0,\alpha_0} $}, which is
also denoted by $\NT v \NT_{0}$.
It is easy to see that the norms
$\NT v \NT_{0}$ and $\NT v \NT$
are
equivalent and satisfy the relation
\begin{equation}
\label{norm_equi}
\underline{{\rm C}}  \NT v \NT^2\;
\leq\;  \NT v \NT_{0}^2 \;\leq \;\overline{{\rm C}}  \NT v \NT^2 ,
\end{equation}
where
\begin{equation}
\overline{{\rm C}} := \max\limits_{ i \in \{ 1,2,3 \} } \frac{\bup{i}}{\blow{i}-\delta}
\quad \textrm{and} \quad
\underline{{\rm C}} := \min\limits_{ i \in \{ 1,2,3 \} } \frac{\blow{i}}{\bup{i}+\delta_i}
\end{equation}
These constants  $\overline{{\rm C}}$ and $\underline{{\rm C}}$
depend only on the problem data and indeterminacy range. They play
an important role in our analysis.

Now, we can define the quantity we are interested in:
\begin{equation}
r_{\indset} :=
\sup\limits_{\tilde u \in \indset} \NT u_0 - \tilde u \NT_{0}.
\end{equation}
A normalized counterpart of $r_\indset$
is defined by the relation
\be
\hat r_{\indset} :=
\sup\limits_{\tilde u \in \indset} \frac{ \NT u_0 - \tilde u \NT_0  }
{  \NT u_0 \NT_0  }.
\ee
\section{Lower bound of $r_{\indset}$}
Problem ${\mathcal P}$ has a variational statement and the solution $u$
can be considered as a minimizer of the functional
\be
J(v):=\frac12 a(v,v)-l(v)
\ee
on the set $V_0$. Using this statement, we can easily derive
 computable lower bounds of the difference between
$u$ and an arbitrary function $v \in V_0$ in terms of the energy norm
(see \cite{RE} where such bounds have been derived for a wide
class of variational problems).
First, we use the identity
\begin{equation}
\label{2.1}
\NT u - v \NT^2 = a(u-v,u-v) = 2 ( J(v) - J(u) ),
\end{equation}
which for quadratic functionals was established in \cite{Mi1}.
Let $w$ be an arbitrary function in $V_0$. Then,
\be
J(v)-J(u) \geq  J(v) - J(v+w)
\ee
and by (\ref{2.1}) we conclude that
\begin{equation}
\label{2.2}
\NT u - v \NT^2  \geq  -a(w+2v,w) + 2 l(w) \quad \forall w \in V_0 .
\end{equation}
We note that for $w = u-v$ the estimate (\ref{2.2})
holds as equality, so that there is no ``gap'' between the left- and
right-hand sides of (\ref{2.2}). The right hand side of (\ref{2.2})
is explicitly computable. It provides the so--called
functional error minorant, which we denote by
 $\mathcal{M}_\ominus^{A,\rho,\alpha}(v,w)$
 (if no confusion may arise, we also
 use a simplified notation $\mathcal{M}_\ominus(v,w)$).
 This functional serves as the main tool when deriving the lower bound
 of $r_{\indset}$.

\begin{Theorem}
\label{Th2.1}
Assume that all the assumptions of Section \ref{intro} hold. Then
\begin{equation}
\label{2.3}
r_\indset^2  \geq  \underline{{\rm C}}
\sup_{w \in V_0} {\mathrm M}^{r_\indset}_\ominus (u_0,w) ,
\end{equation}
where $w$ is an arbitrary function in $V_0$ and
\begin{multline}
\label{2.4}
{\mathrm M}^{r_\indset}_\ominus (u_0,w) :=  -\NT w \NT_{0}^2
+ \delta_1 \IntO  | \nabla w + 2 \nabla u_0 | \; | \nabla w | \; dx + \\
+ \delta_2 \IntO  | ( w + 2 u_0 ) w | \; dx
+ \delta_2 \IntG{3}  | ( w + 2 u_0 ) w | \; ds .
\end{multline}
\end{Theorem}
\proofbegin
We have
\begin{equation}
\label{2.5}
r_{\indset}=
\sup\limits_{\tilde u \in \indset} \NT u_0 - \tilde u \NT_{0}\geq
\underline {C} \sup\limits_{\tilde u \in \indset} \NT u_0 - \tilde u \NT.
\end{equation}
On the other hand
\be
 \sup\limits_{\tilde u \in \indset} \NT u_0 - \tilde u \NT^2 & = &
 \sup\limits_{\tilde u \in \indset}  \left\{ \sup\limits_{w \in V_0}
  \mathcal{M}_\ominus(u_0,w) \right\} \\
& = &
 \sup\limits_{w \in V_0}  \left\{ \sup\limits_{A \in \Lambda_A, \rho \in \Lambda_\rho, \alpha \in \Lambda_\alpha}  \mathcal{M}_\ominus^{A,\rho,\alpha}(u_0,w) \right\} .
\ee
and we conclude that
\begin{equation}
 \label{2.6}
r_{\indset}^2 \geq \underline{{\rm C}}
 \sup\limits_{w \in V_0}  \left\{
\sup\limits_{
A \in \Lambda_A, \rho \in \Lambda_\rho, \alpha \in \Lambda_\alpha
}
\mathcal{M}_\ominus^{A,\rho,\alpha}(u_0,w)
\right\} .
\end{equation}
Now our goal is to estimate the right-hand side of (\ref{2.6})
from below. For this purpose, we exploit the structure of the minorant,
which allows to explicitly
evaluate effects caused by indeterminacy of the coefficients.


It is easy to see that the minorant can be represented
as follows:
\begin{multline}
\label{2.7}
{\mathcal M}_\ominus^{A,\rho,\alpha}(u_0,w) = - \IntO (A_0+\delta_1 \Psi) (\nabla w + 2 \nabla u_0 ) \cdot \nabla w \, dx \\
 - \IntO (\rho_0+\delta_2 \psi_\rho) (w + 2 u_0 ) w \, dx \\
 - \IntG{3} (\alpha_0+\delta_3 \psi_\alpha) (w + 2 u_0 ) w \, ds + 2 l(w)
\end{multline}
Note that
\be
\IntO ( A_0 \nabla u_0 \cdot \nabla w  dx + \rho_0 u_0w ) dx +
\IntG{3} \alpha_0 u_0w ds=l(w).
\ee
Hence,
\begin{multline}
\label{2.8}
{\mathcal M}_\ominus^{A,\rho,\alpha}(u_0,w) = - \IntO A_0 \nabla w \cdot \nabla w \, dx - \IntO \rho_0 w^2 \, dx - \IntG{3} \alpha_0 w^2 \, ds \\
- \delta_1 \IntO  \Psi ( \nabla w + 2 \nabla u_0 ) \cdot  \nabla w  \; dx
- \delta_2 \IntO  \psi_\rho ( w + 2 u_0 ) w  \; dx\\
- \delta_3 \IntG{3} \psi_\alpha  ( w + 2 u_0 ) w  \; ds
\end{multline}
and we obtain
\begin{multline}
\label{2.9}
\sup\limits_{
A \in \Lambda_A, \rho \in \Lambda_\rho, \alpha \in \Lambda_\alpha
}
\mathcal{M}_\ominus^{A,\rho,\alpha}(u_0,w)=-\NT w\NT^2_0+\\
+\delta_1
 \sup\limits_{
|\Psi|\leq 1}
  \IntO  \Psi ( \nabla w + 2 \nabla u_0 ) \cdot  \nabla w  \; dx
+ \delta_2
\sup\limits_{
|\psi_\rho|\leq 1
}
\IntO  \psi_\rho ( w + 2 u_0 ) w  \; dx\\
+ \delta_3
\sup\limits_{
|\psi_\alpha|\leq 1
}
\IntG{3} \psi_\alpha  ( w + 2 u_0 ) w  \; ds.
\end{multline}
 The integrand
of the first integral in the right-hand side of (\ref{2.9}) can be
presented as $\Psi:\tau$,
 where
$$
 \tau = \nabla w \otimes (\nabla w +
2 \nabla u_0)
$$
and $\otimes$ stands for the diad product.
For the first supremum we have
\begin{equation}
\label{2.10}
\sup_{|\Psi| \leq 1} \left\{
\IntO \Psi : \tau \; dx \right\} =
\IntO |\tau|
\; dx \; .
\end{equation}
Analogously, we find that
\ben
\label{2.11}
\sup\limits_{
|\psi_\rho|\leq 1
}
\IntO  \psi_\rho ( w + 2 u_0 ) w  \; dx\leq
\IntO  | (w + 2 u_0 ) w | dx,
\\
\label{2.12}
\sup\limits_{
|\psi_\alpha|\leq 1
}\IntO  \psi_\rho ( w + 2 u_0 ) w  \; dx\leq
\IntG{3} |( w + 2 u_0 ) w| ds.
\een
By (\ref{2.9})--(\ref{2.12}), we arrive at the relation
\begin{multline}
\label{2.13}
\sup\limits_{A,\rho,\alpha}{\mathcal M}_\ominus^{A,\rho,\alpha}(u_0,w)
 = - \NT w\NT^2_0+ \delta_1 \IntO  |( \nabla w + 2 \nabla u_0 )
  \otimes  \nabla w|  \; dx\\
+ \delta_2 \IntO  | ( w + 2 u_0 ) w | \; dx
+ \delta_3 \IntG{3} |( w + 2 u_0 ) w | \; ds .
\end{multline}
which together with (\ref{2.6}) leads to \eqref{2.3}.
\proofend

\vspace{10pt}
Theorem \ref{Th2.1} gives a general form of the lower bound
of $r_\indset$. Also, it creates a basis for practical computation
of this quantity. Indeed, let $V_{0h}\subset V_0$ be a finite dimensional
space. Then
\begin{equation}
\label{2.14}
r_\indset^2  \geq  \underline{{\rm C}}
\sup_{w \in V_{0h}} {\mathrm M}^{r_\indset}_\ominus (u_0,w).
\end{equation}
It is worth noting that the wider set $V_{0h}$ we take the better
lower bound of the radius we compute. However, as it is shown below,
a meaningful lower bound can be deduced even analytically.
\begin{Corollary}
 \label{lb_rough}
  Under assumptions of  Theorem \ref{Th2.1},
  \begin{equation}
  \label{2.15}
    r_\indset^2
    \geq
    \underline{{\rm C}} \; r_{\indset \ominus}^2
    \quad \textrm{ and } \quad
    \hat r_{\indset}^2
    \geq
    \underline{{\rm C}} \; \hat r_{\indset \ominus}^2 ,
  \end{equation}
  where
  \begin{equation}
  \label{2.16}
     r_{\indset \ominus}^2 = \frac{\NT u_0 \NT^4_\delta }
    { \NT u_0\NT_0 ^2 -  \NT u_0\NT_\delta ^2  }\geq  \frac{\Theta^2}{1-\Theta}  \NT u_0 \NT_0 ^2 ,
  \end{equation}
  where
  \be
  \NT u_0 \NT^2_\delta:=
\delta_1 ||\nabla u_0||_\Omega^2 + \delta_2 ||u_0||_\Omega^2 + \delta_3 ||u_0||_{\Gamma_3}^2
  \ee
and
\be
\Theta:=\min\limits_{ i \in \{ 1,2,3 \} } \frac{\delta_i}{\bup{i}}.
\ee
For the normalized radius, we have
  \begin{equation} \label{2.17}
    \hat r_{\indset \ominus}^2 =\frac{\Theta^2}{1-\Theta}.
  \end{equation}
\end{Corollary}

\proofbegin Use \eqref{2.3} and set
\begin{equation}
 \label{2.18}
w = \lambda u_0 ,
\end{equation}
where $\lambda \in {\mathbb R}$. Then we observe that,
\ben
\label{2.19}
&& r_{\indset}^2 \geq \underline{{\rm C}} \Bigg( - \lambda^2 \NT u_0 \NT_0 ^2
+ \lambda (\lambda + 2) \NT u_0 \NT^2_\delta \Bigg) \; .
\een
The right hand side of (\ref{2.19}) is a quadratic function with respect to $\lambda$.
It attains its maximal value if
\be
\lambda \NT u_0 \NT^2_0=(\lambda+1)\NT u_0 \NT^2_\delta,
\ee
i.e., if for
\be
\lambda=\frac{\NT u_0 \NT^2_\delta}{\NT u_0 \NT^2_0-\NT u_0 \NT^2_\delta} .
\ee
 Substituting this $\lambda$, we arrive at \eqref{2.16}. Note that
\begin{multline}
 \NT u_0 \NT_0 ^2 = \IntO ( A_0 \nabla u_0 \cdot \nabla u_0 + \rho_0 u_0^2\; dx + \IntG{3} \alpha_0 u_0 \; ds
\geq \\
\geq \IntO ( \blow{1} \nabla u_0 \cdot \nabla u_0 + \blow{2} u_0^2 ) \,dx +
\IntG{3} \blow{3} u_0 \,ds
>  \\
>\delta_1 ||\nabla u_0||_\Omega^2 + \delta_2 ||u_0||_\Omega^2 + \delta_3 ||u_0||_{\Gamma_3}^2=
\NT u_0 \NT^2_\delta,
\end{multline}
so that $\lambda$ (and the respective lower bound) is positive.
Moreover,
\begin{multline}
\label{2.21}
\NT u_0 \NT^2_\delta
\geq
\frac{\delta_1}{\bup{1}} \IntO A_0 \nabla u_0 \cdot \nabla u_0 \; dx
+ \frac{\delta_2}{\bup{2}} \IntO  \rho_0 u_0^2 \; dx
+ \frac{\delta_3}{\bup{3}} \IntG{3}  \alpha_0 u_0^2 \; ds \\
\geq
\Theta \NT u_0 \NT_0 ^2.
\end{multline}
Also,
\begin{multline}
\label{2.22}
\NT u_0\NT_0 ^2 -  \NT u_0\NT_\delta ^2\\=
\IntO (A_0-\delta_1I) \nabla v \cdot \nabla v \; dx + \IntO (\rho_0-\delta_2) v^2 dx +
\IntG{3} (\alpha_0-\delta_3) v^2 \, ds\\
\geq
\left(1-\frac{\delta_1}{\bup{1}}\right)
\IntO A_0\nabla v \cdot \nabla v \, dx + \left(1-\frac{\delta_2}{\bup{2}}\right)
\IntO \rho_0 v^2 dx +
\left(1-\frac{\delta_3}{\bup{3}}\right)\IntG{3} \alpha_0 v^2 \, ds\\
\geq \max\limits_{i=1,2,3}
\left(1-\frac{\delta_i}{\bup{i}}\right)\NT u_0\NT^2_0=(1-\Theta)\NT u_0\NT^2_0.
\end{multline}
By (\ref{2.21}) and (\ref{2.22}), we arrive at relation
\begin{equation}
\label{2.23}
r_{\indset \ominus}^2 \geq  \frac{\Theta^2}{1-\Theta}  \NT u_0 \NT_0 ^2,
\end{equation}
which implies (\ref{2.16}) and (\ref{2.17}).

\section{Upper bound of $r_{\indset}$}
A computable upper bound
of $r_{\indset}$ can be derived with the help of a posterioiri error majorant
of the functional type, which are derived by purely functional methods without
attracting any information on the mesh and method used. For a wide class of problems they were
derived in \cite{RE,RE1,RE2} by variational techniques and in \cite{RE2,RS,RST}
by transformations of integral identities (see also the papers cited therein).
Below
we  derive functional error majorant for our class of problems using the latter method
based on transformations  of the respective integral identity.
After that, we use it's properties to derive the desired upper bound in \ref{der_ub}.

\subsection{Error majorant}
\label{der_maj}
Let $v\in V_0$ be an admissible approximation of the exact solution $u$
(generated by $A$, $\varrho$, and $\alpha$).
From \eqref{varform} it follows that for any $w\in V_0$
\begin{multline}
a(u-v,w)= \IntO fw \; dx + \IntG{2} Fw \; ds + \IntG{3} Gw \; ds - \\
    -\IntO A \nabla v \cdot \nabla w \; dx - \IntO \rho v w \; dx - \IntG{3} \alpha vw \; ds + \\
+ \IntO \left( \dvg(y) w + y \cdot \nabla w \right) \; dx - \int\limits_{\Gamma_2\cup\Gamma_3}
( y\cdot\nu)\, w \; ds\,,
\end{multline}
where $\nu$ denotes unit outward normal to $\Gamma$
and
\be
y\in H^+(\dvg,\Omega):=
\{y\in H(\dvg,\Omega)\,\mid\,y\cdot\nu\in L^2(\Gamma_2\cup\Gamma_3)\}.
\ee
We note that the last line is zero for all $y \in H(\Omega,\dvg)$
(in view of the integration-by-parts formula). We regroup the terms and rewrite the relation
as follows:
\begin{equation}
a(u-v,w)= I_1 + I_2 + I_3 + I_4 ,
\end{equation}
where
\be
I_1 & := & \IntO r_1(v,y)w \; dx := \IntO (f - \rho v + \dvg(y) ) w \; dx ,\\
I_2 & := & \IntO r_2(v,y)w \; dx := \IntG{3} (G-\alpha v-y\cdot\nu\,) w \; ds ,\\
I_3 & := & \IntG{2} (F-y\cdot\nu\,) w \; ds ,\\
I_4 & := & \IntO (y - A \nabla v) \cdot \nabla w \; dx .
\ee
Now we can estimate each term separately by the Friedrichs and  trace inequalities
(which holds due to our assumptions concerning $\Omega$). We have
\be
|| w ||_\Omega^2 \leq C_1(\Omega) || \nabla w ||_\Omega^2 \quad \forall w \in V_0, \\
|| w ||_{2,\Gamma_2}^2 \leq C_2(\Omega,\Gamma_2) || \nabla w ||_\Omega^2 \quad \forall w \in V_0, \\
|| w ||_{2,\Gamma_3}^2 \leq C_3(\Omega,\Gamma_3) || \nabla w
||_\Omega^2 \quad \forall w \in V_0. \ee When estimating the
integrands of $I_1$ and $I_2$, we  introduce additional functions
$\mu_1$ and $\mu_2$, which have values in  $[0, 1]$.
\begin{multline*}
I_1  = \IntO \frac{\mu_1}{\sqrt{\rho}} r_1(v,y) \sqrt{\rho} w \; dx + \IntO (1-\mu_1) r_1(v,y) w \; dx \\
     \leq  || \frac{\mu_1}{\sqrt{\rho}} r_1(v,y) ||_\Omega \; || \sqrt{\rho} w ||_\Omega +
             \sigma_1||  (1-\mu_1) r_1(v,y) ||_\Omega \;  \left(\IntO A\nabla w\cdot\nabla w \,dx\right)^{1/2}
\end{multline*}
and
\be
I_2 & = & \IntG{3} \frac{\mu_2}{\sqrt{\alpha}} r_2(v,y) \sqrt{\alpha} w \; dx + \IntG{3} (1-\mu_2) r_2(v,y) w \; dx \\
    & \leq & || \frac{\mu_2}{\sqrt{\alpha}} r_2(v,y) ||_{\Gamma_3} \;
     || \sqrt{\alpha} w ||_{\Gamma_3} +\\
&&+             ||  (1-\mu_2) r_2(v,y) ||_{\Gamma_3} \;  \sigma_3
             \left(\IntO A\nabla w\cdot\nabla w \,dx\right)^{1/2},
\ee
and
\begin{equation}
I_3 \leq || F-y\cdot\nu\, ||_{\Gamma_2}
\sigma_2
\left(\IntO A\nabla w\cdot\nabla w \,dx\right)^{1/2},
\end{equation}
where
\be
\sigma_1=\sqrt{ \frac{C_1(\Omega)}{\blow{1}} },\quad
\sigma_2=\sqrt{ \frac{C_1(\Omega)C_2(\Omega,\Gamma_2)}{\blow{1}} },
{\rm and}\quad
\sigma_3=\sqrt{ \frac{C_1(\Omega)C_3(\Omega,\Gamma_3)}{\blow{1}} }.
\ee
We note that a similar approach
was used in \cite{RS} for the reaction-diffusion problem and in \cite{RST}
for the generalized Stokes problem.
In these publications it was shown that splitting of the
residual term (performed with the help of a single function $\mu$) allows to obtain error majorants
that are insensitive with respect to small values of the lower term coefficient and at the same
time sharp (i.e., have no irremovable gap between the left- and right-hand sides). In our case,
we have two lower terms, so that we need two functions $\mu_1$ and $\mu_2$ to split the respective
residual terms.

The term $I_4$ is estimated as follows:
\begin{equation}
I_4 \leq D(\nabla v,y)^{\frac{1}{2}} \left(\IntO A\nabla w\cdot\nabla w \,dx\right)^{1/2},
\end{equation}
where
\begin{equation}
D(\nabla v, y) = \IntO ( y-A\nabla v ) \cdot ( \nabla v - A^{-1} y ) \; dx .
\end{equation}
We collect all the terms and obtain
\begin{multline}
a(u-v,w) \leq \Big( D(\nabla v,y)^{1/2} +
\sigma_1 ||  (1-\mu_1) r_1(v,y) ||_\Omega + \\
+ \sigma_3 || (1-\mu_2) r_2(v,y) ||_{\Gamma_3} + \sigma_2 || F-y\cdot\nu\, ||_{\Gamma_2}  \Big) \left(\IntO A\nabla w\cdot\nabla w \,dx\right)^{1/2} + \\
 + || \frac{\mu_1}{\sqrt{\rho}} r_1(v,y) ||_\Omega \; || \sqrt{\rho} w ||_\Omega +
 || \frac{\mu_2}{\sqrt{\alpha}} r_2(v,y) ||_{\Gamma_3} \; || \sqrt{\alpha} w ||_{\Gamma_3} .
\end{multline}
Set $w=u-v$ and
use the Cauchy-Schwartz inequality
\begin{equation}
\label{CS}
\sum\limits_{i=1}^d x_i y_i  \leq \sqrt{ \sum\limits_{i=1}^d  x_i^2}\;
 \sqrt{ \sum\limits_{i=1}^d {y_i^2}}.
\end{equation}
Then, we arrive at the estimate
\begin{multline}
\label{3.8}
\NT u-v\NT ^2 \leq \Big( D(\nabla v,y)^{1/2} + \sigma_1 ||  (1-\mu_1) r_1(v,y) ||_\Omega + \\
+ \sigma_3 || (1-\mu_2) r_2(v,y) ||_{\Gamma_3} + \sigma_2 || F-y\cdot\nu\, ||_{\Gamma_2}  \Big)^2 \\
 + || \frac{\mu_1}{\sqrt{\rho}} r_1(v,y) ||_\Omega^2 + || \frac{\mu_2}{\sqrt{\alpha}} r_2(v,y) ||_{\Gamma_3}^2  .
\end{multline}
It is worth remarking, that the estimate (\ref{3.8}) provides a guaranteed upper bound
of the error for {\em any} conforming approximation of the problem (\ref{basic1}--\ref{basic4}). The
estimate has a form typical for all functional a posteriori estimates: it is presented by the sum
of residuals of the basic relations with multipliers that depend on the constants in the respective
functional  (embedding) inequalities for the domain and boundary parts.

However, for our subsequent goals, it is desirable to have the majorant in a form that
involve only quadratic terms. Such a form can be easily derived from (\ref{3.8})
if we square both parts and apply the algebraic inequality
  \eqref{CS} to the first term (with multipliers $\gamma_i>0$, $i=1,2,3,4$). Then, we obtain
\begin{multline}
\label{3.9}
\NT u-v\NT ^2 \leq \kappa
 \Bigl( \gamma_1 D(\nabla v,y) + \gamma_2 ||(1-\mu_1) r_1(v,y) ||_\Omega^2\\
+ \gamma_3 || (1-\mu_2) r_2(v,y) ||_{\Gamma_3}^2 + \gamma_4 || F-y\cdot\nu\, ||_{\Gamma_2}^2  \Bigr) + \\
 + || \frac{\mu_1}{\sqrt{\rho}} r_1(v,y) ||_\Omega^2 + || \frac{\mu_2}{\sqrt{\alpha}} r_2(v,y) ||_{\Gamma_3}^2  .
\end{multline}
where
\be
\kappa :=\frac{1}{\gamma_1} + \frac{\sigma^2_1}{ \gamma_2} +
\frac{\sigma^2_2}{ \gamma_3} +  \frac{\sigma^2_3}{\gamma_4}.
\ee
We note that (\ref{3.9}) coincides with (\ref{3.8}) if
\ben
\gamma_1=\bar\gamma_1 & := & D(\nabla v,y)^{-1/2} , \\
\gamma_2=\bar\gamma_2  & := & \frac{  \sigma_1  }{ ||(1-\mu_1) r_1(v,y) ||_\Omega } , \\
\gamma_3=\bar\gamma_3  & := & \frac{ \sigma_2 }{ ||(1-\mu_2) r_2(v,y) ||_{\Gamma_3} } , \\
\gamma_4=\bar\gamma_4 &: = &  \frac{ \sigma_3 }{ || F-y\cdot\nu\, ||_{\Gamma_2} } .
\een
Certainly, the estimate (\ref{3.9}) looks more complicated with respect
to (\ref{3.8}). However, it has an important advantage: the
weight functions $\mu_1$ and $\mu_2$ enter it as quadratic integrands, so that we can
easily find their optimal form adapted to a particular $v$ and the respective
error distribution.

In the simplest case, we take $\mu_1=\mu_2=0$, which yields the estimate
\begin{multline} \label{3.14}
\NT u-v\NT ^2 \leq \kappa
\times \Big( \gamma_1 D(\nabla v,y) + \gamma_2 || r_1(v,y) ||_\Omega^2 +\\
+\gamma_3 || r_2(v,y) ||_{\Gamma_3}^2 + \gamma_4 || F-y\cdot\nu\, ||_{\Gamma_2}^2  \Big).
\end{multline}
Another estimate arises if we set  $\mu_1=\mu_2=1$.
In this case, the terms with factors $\sigma_1$ and $\sigma_3$ in (\ref{3.8})
are equal to zero, so that subsequent relations do not contain the terms with
multipliers formed by $\gamma_2$ and $\gamma_3$. Hence, we arrive at the estimate
\begin{multline}
\label{3.15}
\NT u-v\NT ^2 \leq \Big( \frac{1}{\gamma_1} +  \frac{\sigma^2_3}{\gamma_4} \Big) \times \Big( \gamma_1 D(\nabla v,y) + \gamma_4 || F-y\cdot\nu\, ||_{\Gamma_2}^2  \Big) + \\
 + || \frac{1}{\sqrt{\rho}} r_1(v,y) ||_\Omega^2 + || \frac{1}{\sqrt{\alpha}} r_2(v,y) ||_{\Gamma_3}^2  .
\end{multline}
The estimate \eqref{3.15} involves ``free'' parameters $\gamma_1$ and $\gamma_4$
and a ´´free'' vector--valued function $y$ (which can be thought of as an image of the true flux).
There exist a combination of these free parameters that makes the left-hand side
of the estimate equal to the right-hand one. Indeed, set $y=A \nabla u$. Then
\be
r_1(v,y) & = & \rho (u-v) \quad{\rm in}\;\Omega , \\
r_2(v,y) & = & \alpha (u-v) \quad{\rm on}\;\Gamma_3 , \\
F - y\cdot\nu & = & 0\quad{\rm on}\;\Gamma_2
\ee
and we find that for $\gamma_4$ tending to infinity and for any $\gamma_1>0$ the right-hand side of (\ref{3.15})
coincides with the energy norm of the error.
However, the estimate (\ref{3.15}) has a drawback: it is sensitive
with respect to $\rho$ and $\alpha$
and may essentially overestimate the error if $\rho$ or $\alpha$ are small.
On the other hand, the right-hand side of \eqref{3.14} is stable with respect to small values of
$\rho$ and $\alpha$. Regrettably, it does not possess the ``exactness'' (in the above-discussed sense)
because it may have a ``gap'' between the left- and right-hand sides for any $y$.

An upper bound of the error that combines positive features of (\ref{3.14}) and (\ref{3.15})
can be derived (as in \cite{RS,RST}) if a certain optimization procedure is used in order to
select optimal functions $\mu_1$ and $\mu_2$.
If $y$ is given, then optimal $\mu_1$ and $\mu_2$ can be found analytically. It is not difficult
to see that  $\mu_1$ must minimize the quantity
\begin{equation}
\IntO\left(
\kappa \gamma_2 (1-\mu_1)^2 + \frac{\mu_1^2}{\rho}  \right) r_1(v,y)^2 dx.
\end{equation}

The quantity attains its minimum with
\begin{equation}
\mu_1(x)=\mu_1^{opt}(x): = \frac{ \kappa \gamma_2 }{ \kappa \gamma_2 + \rho(x)^{-1}}  \quad \textrm{ in } \Omega .
\end{equation}
Similarly, we find that the integrals associated with $\Gamma_3$ attain minimum if
\begin{equation}
\mu_2(x)=\mu_2^{opt}(x) := \frac{ \kappa \gamma_3}{ \kappa \gamma_3+ \alpha(x)^{-1} } \quad \textrm{ on } \Gamma_3 .
\end{equation}
Substituting these values to \eqref{3.9} results in the estimate
\begin{multline}
\label{3.19}
\NT u-v\NT ^2 \leq \kappa  \Big(  \gamma_1 D(\nabla v , y)
+ \gamma_2 || \frac{ \sqrt{ \kappa^2 \gamma_2^2 \rho + 1 } }{ \kappa \gamma_2 \rho   + 1 }  \; r_1(v,y) ||_\Omega^2 + \\
+ \gamma_3 ||  \frac{ \sqrt{ \kappa^2 \gamma_3^2 \alpha + 1 } }{ \kappa \gamma_3 \alpha   + 1 } \; r_2(v,y) ||_{\Gamma_3}^2
 + \gamma_4 || F-y\cdot\nu\, ||_{\Gamma_2}^2  \Big) .
\end{multline}

\begin{Remark}
For practical computations, it may be easier to use \eqref{3.9} and directly
minimize its right-hand side with respect to $\mu_i$, $\gamma_i$, and $y$ using the following
iteration procedure:
\begin{enumerate}
\item Keep $\gamma_i$ and $\mu_j$ fixed in \eqref{3.9} and minimize resulting quadratic functional of $y$ in suitable finite subspace. This task can be reduced to solving a system of linear equations.
\item  Compute $\gamma_i^{opt}$.
\item Compute $\mu_j^{opt}$ and repeat from step 1.
\end{enumerate}
\end{Remark}
We denote the right hand side of \eqref{3.9}  by ${\mathcal{M}_\oplus(v,y,\gamma,\mu_1,\mu_2)}$.
This error majorant provides a guaranteed upper bound of the error,
i.e.,
\begin{equation}
\label{3.20}
\NT u-v\NT ^2 \leq {\mathcal{M}_\oplus(v,y,\gamma,\mu_1,\mu_2)}.
\end{equation}
It is exact (in the above-discussed sense). Indeed, for
$y=A\nabla u$ and $\mu_1=\mu_2=1$ we obtain
\begin{equation}
\label{3.21}
\inf\limits_{\gamma_i>0}{\mathcal{M}_\oplus(v,A \nabla u,\gamma,1,1)}= \NT u-v\NT ^2 .
\end{equation}
Also, it directly follows from the structure of (\ref{3.19}) that the right-hand side
is insensitive to small values of $\rho$ and $\alpha$.
\subsection{Upper bound}
\label{der_ub}
In this section, we
derive an upper bound of $r_{\indset}$. For this purpose, we  use the majorant
${\mathcal{M}_\oplus(v,y,\gamma,\mu_1,\mu_2)}$. Since the majorant
nonlinearly depends on $A$, $\rho$, and $\alpha$, taking supremum over the respective
indeterminacy sets imposes a more complicated task than that for the minorant.
For a class of diffusion problems this task was solved in \cite{RE,RE3}. Below, we
deduce a simpler estimate, which can be easily exploited in practical computations
and serves as a natural counterpart for the lower bound derived in Corollary 2.1.

\begin{Proposition} \label{comp}
Assume that all the assumptions of Section \ref{intro} hold. Then
\begin{equation} \label{ub00}
r_{\indset}^2 \leq  \overline{{\rm C}} \; r_{\indset \oplus}^2
\quad \textrm{and} \quad
\hat r_\indset^2 \leq \overline{{\rm C}} \; \hat r_{\indset \oplus}^2 ,
\end{equation}
where
\begin{equation}
 \label{3.23}
 r_{\indset \oplus}^2 =  \frac{\delta_1^2}{\blow{1}-\delta_1} ||\nabla u_0 ||^2_\Omega + \frac{\delta_2^2}{\blow{2}-\delta_2} || u_0 ||_\Omega^2 + \frac{\delta_3^2}{\blow{3}-\delta_3} || u_0 ||_{\Gamma_3}^2
\end{equation}
and
\begin{equation}
\label{3.24}
\hat r_{\indset \oplus}^2 =
\max\limits_{i \in \{ 1,2,3\}} \frac{\delta_i^2}{\blow{i} (\blow{i}-\delta_i)}.
\end{equation}
\end{Proposition}
\proofbegin
By properties of the majorant, we have
\be
\sup\limits_{\tilde u \in \indset} \NT u_0-\tilde u\NT^2
& = &
\sup\limits_{\tilde u \in \indset} \left\{ \inf\limits_{y,\mu_i,\gamma_j}  {\mathcal{M}_\oplus^{A,\rho,\alpha}(u_0,y,\gamma,\mu_1,\mu_2)} \right\} \\
& \leq &
\inf\limits_{y,\mu_i,\gamma_j} \left\{ \sup\limits_{A,\rho,\alpha}  {\mathcal{M}_\oplus^{A,\rho,\alpha}(u_0,y,\gamma,\mu_1,\mu_2)} \right\} .
\ee
Applying \eqref{norm_equi}, we obtain
\begin{equation}
\label{3.25}
r_\indset^2 \leq \overline{{\rm C}} \inf\limits_{y,\mu_i,\gamma_j} \left\{ \sup\limits_{A,\rho,\alpha}  {\mathcal{M}_\oplus^{A,\rho,\alpha}(u_0,y,\gamma,\mu_1,\mu_2)} \right\} .
\end{equation}
Our task is to explicitly estimate the  term in brackets.
For this purpose, we estimate from above the last two terms of the majorant and represent
it in the form
\begin{multline} \label{3.26}
{\mathcal{M}_\oplus^{A,\rho,\alpha}(u_0,y,\gamma,\mu_1,\mu_2)} \leq\\ \kappa
\Bigg(  \gamma_1 D(\nabla v , y)
+ \left\| \sqrt{ \gamma_2 \kappa (1-\mu_1)^2 +
\frac{\mu_1^2}{\kappa(\blow{2} -\delta_2)}} \; r_1(v,y) \right\|_\Omega^2 + \\
+ \left\|
\sqrt{ \gamma_3 \kappa (1-\mu_2) +
\frac{\mu_2}{\kappa(\blow{3} -\delta_3)}} \;  r_2(v,y) \right\|_{\Gamma_3}^2
 + \gamma_4 || F-y\cdot\nu\, ||_{\Gamma_2}^2  \Bigg) .
\end{multline}
Now, we find upper bounds with respect to
 $A \in \Lambda_A$, $\rho \in \Lambda_\rho$, and $\alpha \in \Lambda_\alpha$ separately.

First, we consider the term $D$ generated by $A$ and $A^{-1}$:
\begin{multline}
\label{3.27}
\sup\limits_{A \in \Lambda_A} D(\nabla u_0, y) =
\sup\limits_{|\Psi|<1} \IntO (A_0+\delta_1 \Psi)^{-1} |(A_0+\delta \Psi)\nabla u_0 -y|^2 \; dx \\
\leq
\frac{1}{\blow{1}-\delta_1} \sup\limits_{|\Psi|<1} \left\{ ||A_0 \nabla u_0 -y ||^2 +
2 \delta_1 \IntO \Psi \nabla u_0 \cdot (A_0 \nabla u_0 - y) \; dx +
\delta_1^2 || \Psi \nabla u_0 ||^2 \right\} \\
\leq
\frac{1}{\blow{1}-\delta_1} \left( ||A_0 \nabla u_0 -y ||^2_\Omega +
2 \delta_1 \IntO |\nabla u_0 | \; | A_0 \nabla u_0 - y | \; dx +
\delta_1^2 || \nabla u_0 ||^2_\Omega \right).
\end{multline}
For the term related to the error in equilibrium equation, we have
\begin{multline}
\label{3.28}
\sup\limits_{\rho \in \Lambda_\rho} || r_1^{\rho}(u_0,y) ||_\Omega^2 =
\sup\limits_{|\psi_2|<1}  \IntO (f- (\rho_0+\delta_2 \psi_2) u_0 + \dvg y )^2 \; dx \\
=
 \sup\limits_{|\psi_2|<1} \IntO \left( \dvg y - \dvg( A_0 \nabla u_0) -\delta_2 \psi_2 u_0 \right)^2 \; dx
\\
\leq
|| \dvg(y-A_0 \nabla u_0) ||_\Omega^2 + 2 \delta_2 \IntO | \dvg(y-A_0 \nabla u_0) | \; |u_0| \; dx
 + \delta_2 || u_0 ||^2.
\end{multline}
Similarly, for  the term related to the error in Robin boundary condition we have
\begin{multline}
\label{3.29}
\sup\limits_{\alpha \in \Lambda_\alpha} ||  r_2^{\alpha}(u_0,y) ||_{\Gamma_3}^2 \leq
 \left\| \frac{\partial (y-A_0 \nabla u_0)}{\partial \nu} \right\|_{\Gamma_3}^2 \\
 + 2 \delta_3 \IntG{3}\left| \frac{\partial (y-A_0 \nabla u_0)}{\partial \nu} \right| \; | u_0 | \; ds + \delta_3^2 ||u_0||_{\Gamma_3}^2.
\end{multline}
It is clear, that  for $y=y_0:=A_0 \nabla u_0$, the estimates (\ref{3.27})--(\ref{3.29})
 attain minimal values. In addition, we
set  in \eqref{3.26} $\mu_1=\mu_2 =1$ and find that
\begin{multline} \label{3.30}
{\mathcal{M}_\oplus^{A,\rho,\alpha}(u_0,A_0\nabla u_0,\gamma,1,1)}\leq\\
\leq
\kappa\left(
 \frac{\delta_1^2 \gamma_1}{\blow{1}-\delta_1}  || \nabla u_0 ||_\Omega^2
+ \frac{\delta_2^2}{\blow{2} -\delta_2}  || u_0 ||_\Omega^2 +
 \frac{\delta_3^2}{\blow{3} -\delta_3}  || u_0 ||_{\Gamma_3}^2\right) .
\end{multline}
Now we tend $\gamma_2$,$\gamma_3$ and $\gamma_4$ (which are contained in $\kappa$) to infinity.
Then, \eqref{3.30} and \eqref{3.25} imply \eqref{3.23}.
An upper bound for the normalized radius follows from the relation
\begin{multline*}
{\mathcal{M}_\oplus^{A,\rho,\alpha}(u_0,A_0 \nabla u_0,\gamma,1,1)} \leq \\
\leq  \frac{\delta_1^2}{\blow{1} (\blow{1}-\delta_1)}
\IntO A \nabla u_0 \cdot \nabla u_0 \; dx +
\frac{\delta_2^2}{\blow{2} (\blow{2}-\delta_2)} || \sqrt{\rho_0} u_0 ||_\Omega^2 +
\frac{\delta_3^2}{\blow{3} (\blow{3}-\delta_3)} || \sqrt{\alpha_0} u_0 ||_{\Gamma_3}^2 \leq \\
\leq \max\limits_{i \in \{ 1,2,3\}} \frac{\delta_i^2}{\blow{i} (\blow{i}-\delta_i)} \NT  u_0 \NT ^2,
\end{multline*}
which leads to \eqref{3.24}. \proofend


\begin{Remark}
\label{b_ok}
The normalized lower bound in Corollary 2.1 is
less than the upper bound established in Proposition 3.1. Indeed, using
\begin{equation*}
 \left(
   1 - \min\limits_i \frac{\delta_i}{\bup{i}}
 \right)^{-1}
=
 \left(
  \max\limits_i \left( 1 - \frac{\delta_i}{\bup{i}} \right)
 \right)^{-1}
=
 \left(
  \min\limits_i \frac{\bup{i}}{\bup{i} - \delta_i}
 \right)
\end{equation*}
to $\hat r_{\indset \ominus}^2$, we arrive at the following relation between bounds:
\begin{multline*}
\frac{  \overline{{\rm C}} \; \hat r_{\indset \oplus}^2 }{  \underline{{\rm C}} \; \hat r_{\indset \ominus}^2  } =
\frac{ \left( \max\limits_{i} \frac{\bup{i}}{\blow{i} - \delta_i} \right) \left(  \max\limits_{i} \frac{\delta_i^2}{\blow{i} (\blow{i}-\delta_i)} \right) }{ \left(  \min\limits_{i} \frac{\blow{i}}{\bup{i} + \delta_i} \right) \left( \min\limits_{i} \frac{\delta_i^2}{\bup{i}^2} \right) \left(  \min\limits_{i} \frac{\bup{i}}{{\bup{i}-\delta_i}} \right)  } = \\
= \left( \max\limits_{i } \frac{\bup{i}}{\blow{i} - \delta_i} \right) \left(  \max\limits_{i} \frac{\delta_i^2}{\blow{i} (\blow{i}-\delta_i)} \right)  \left(  \max\limits_{i } \frac{\bup{i} + \delta_i}{\bup{i}} \right) \left( \max\limits_{i} \frac{\bup{i}^2}{\delta_i^2} \right) \left(  \max\limits_{i} \frac{{\bup{i}-\delta_i}}{\bup{i}} \right) .
\end{multline*}
Maximums can be estimated from below by setting $i=1$ everywhere. The expression simplifies to
\begin{equation}
  \frac{  \overline{{\rm C}} \; \hat r_{\indset \oplus}^2 }{  \underline{{\rm C}} \; \hat r_{\indset \ominus}^2  }
  \geq
  \frac{\bup{1}}{\blow{1}}
  \left(
    \frac{\bup{1}+\delta_1}{\blow{1}-\delta_1}
  \right) \geq 1 .
\end{equation}
\end{Remark}


\end{document}